\documentclass[12pt,a4paper,reqno]{amsart}
\usepackage[centertags]{amsmath}
\usepackage{amsfonts}
\usepackage{amssymb}
\usepackage{amsthm}
\usepackage{a4}
\usepackage{enumerate}

\usepackage{bbm}

\newlength{\defbaselineskip}
\newcommand{\setlinespacing}[1]%
           {\setlength{\baselineskip}{#1 \defbaselineskip}}

\begin{document}

\title[Non-point symmetry  reduction method of partial ... ]{Non-point symmetry  reduction method of partial differential equations }

\author[I. Tsyfra]{Ivan Tsyfra$^{1}$}
\author[T. Czy\.zycki]{Tomasz Czy\.zycki$^{2}$}

%\date{\today}

\maketitle
\noindent $^1$ AGH University of Science and Technology, Faculty of Applied Mathematics, 30 Mickiewicza Avenue, 30-059 Krakow, Poland \\[4pt]
 $^2$ Institute of Mathematics, University of Bia\l ystok,
Cio\l kowskiego 1M, 15-245 Bia\l ystok, Poland \\[4pt]
Email: tsyfra@agh.edu.pl, tomczyz@math.uwb.edu.pl \\

\noindent MSC classification numbers: 35A20, 35Q55, 35Q72, 82D75  \\[4pt]

\begin{abstract}
We  study the symmetry reduction of nonlinear partial differential equations with two independent variables. We propose new ans\"atze reducing nonlinear evolution equations to system of ordinary differential equations. The ans\"atze are constructed
by using operators of non-point classical and conditional symmetry. Then we find solution to nonlinear heat equation which can not be obtained in the framework of the classical Lie approach. By using operators of Lie--B\"acklund symmetries we construct the solutions of nonlinear hyperbolic equations depending on  arbitrary smooth function of one variable too. We show that the method can be applied to nonevolutionary partial 
differential equations.

\end{abstract}

\vspace{0,2cm}

\section{Introduction}
It is well known that the classical Lie symmetry method of  point transformations is often used for reducing the number
of independent variables in partial differential equation to obtain ordinary differential equations. After integration of reduced differential equations one can obtain partial solutions of the equation under study
\cite{BK,O,Ovs}. The main problem is that the maximal invariance group of point transformations of differential equations used in applications is not sufficiently wide and thus the group approach can not be successfully applied to these equations. The concept of generalized conditional symmetry has been introduced in \cite{Fok, Zhd} to extend the applicability of the symmetry method to the construction of solutions of evolution equations.  The relationship of generalized conditional symmetries of evolution equations to compatibility of system of differential equations is studied in \cite{Kun}.
The method for construction nonlocally related partial differential equation systems for a given partial differential equation has been proposed in \cite{BluKyt}. The starting point for the method is the existence of operator of point symmetry of the equation under study. Through a nonlocally related systems one can construct operators of nonlocal symmetry and nonlocal conservation laws of initial equation.   

 We use operators of non-point classical and conditional symmetries to extend the class of differential equations to which  the symmetry method is applicable. In this paper we study the symmetry reduction of partial differential equations with two independent variables  by using the operators of non--point symmetry because the prolongated operators of classical point symmetry lead to the classical invariant solutions. The method can be naturally generalized  to the multidimensional case. We construct the  ansatz for dependent variable $u$ or its derivatives which reduces the scalar partial differential equation to a system of ordinary differential equations. We use the operators of the classical point symmetry \cite{BK,O} of the corresponding system which are not the prolongated operators of point  symmetries  admitted by the original equation to construct the ansatz for derivatives. The ansatz for $u$ we construct by using ordinary differential equation admitting the operators of Lie--B\"acklund symmetry (in the classical sense \cite{O, Fok}).   We consider nonlinear evolution and wave type equations and present  the operator of conditional symmetry for the corresponding system which generates the B\"acklund transformations for nonlinear wave equation.

  Recall, that the well-known integrable nonlinear differential equations such as Korteweg-de-Vries, sine-Gordon, cubic Schr\"{o}dinger equations admit an infinite number of Lie--B\"acklund symmetry operators \cite{BK,O}. Another goal of this paper is to show  that such important properties of nonlinear partial differential equations as existence of B\"acklund transformations,  linearization, existence of the class of solutions depending on arbitrary function can be related  to their invariance under the finite number of non-point symmetry operators.

\section{Non-point symmetry and reduction of nonlinear wave type and evolution equations with two independent variables }
The concept of differential invariant solutions based on infinite Lie group $G$ is introduced in \cite{Ovs}. This group is a classical symmetry group of point transformations of dependent and independent variables for the equation under study. Generally speaking, analysis similar to that in constructing differential invariant solutions enables us to obtain the ans\"atze for derivatives $u_{x_1}$, $u_{x_2}$ by virtue of  operators of non-point symmetry \cite{Ovs,T2}. 
Let us consider nonlinear differential equation
\begin{equation}\label{2}
u_{x_2x_2}=\frac{1}{1-u_{x_1}^r}, \quad  r \not =0,\pm 1.
\end{equation}

We search for the ansatz for the derivatives of such  form
\begin{equation}\label{1}
\frac{\partial u}{\partial x_1}=R_1(x_1, x_2, u,\varphi_1(\omega) ,\varphi_2(\omega)), \quad \frac{\partial u}{\partial x_2}=R_2(x_1, x_2, u, \varphi_1(\omega), \varphi_2(\omega)) ,
\end{equation}
where $\omega=\omega(x_1, x_2, u)$. Operators of classical and conditional symmetry of the corresponding system can be used to find $R_1$, $R_2$. The corresponding system has the form
\begin{equation}\label{3}
 v^1_2=v^2_1, \quad v^2_2=\frac{1}{1-(v^1)^r},
\end{equation}
where $v^1=u_{x_1}$,$v^2=u_{x_2}$, $v^i_k=v^i_{x_k}$, $i, k=1,2$. To construct ansatz of type (\ref{1}) we use the symmetry operator
\[
 Q=(r+1)x_1\partial_{x_1}+rv^2\partial_{x_2}-v^1\partial_{v^1}+rv^2\partial_{v^2}.
\]
of system (\ref{3}). It is obvious that operator $Q$ generates non-point group transformations for variables $x_1, x_2, u$. It is easy to find the invariants of one-parameter Lie group with generator $Q$
\[
\omega=x_2-v^2, \quad \omega_1=v^1(x_1)^{\frac{1}{r+1}}, \quad \omega_2=v^2(x_1)^{\frac{-r}{r+1}}.
\]
By using these invariants  one can construct the ansatz for $v^1$, $v^2$
\begin{equation}\label{4d}
v^1=(x_1)^{\frac{-1}{r+1}}\varphi_1(\omega), \quad v^2=(x_1)^{\frac{r}{r+1}}\varphi_2(\omega).
\end{equation}
From (\ref{4d}) we have
\begin{equation}\label{5d}
v^2_2=\frac{(x_1)^{\frac{r}{r+1}}\varphi'_2}{1+(x_1)^{\frac{r}{r+1}}\varphi'_2}
\end{equation}
where $\varphi'_2=\frac{d\varphi_2}{d\omega}$.
Substituting (\ref{4d}) and (\ref{5d}) into the equation
\[
v^2_2=\frac{1}{1-(v^1)^r}
\]
yields
\[
(x_1)^{\frac{r}{r+1}}\varphi'_2-\varphi'_2\varphi^r_1=1+(x_1)^{\frac{r}{r+1}}\varphi'_2.
\]
Thus we get the first reduced ordinary differential equation
\begin{equation}\label{6d}
\varphi'_2\varphi^r_1=-1.
\end{equation}
The second one we obtain from the compatibility condition $v^1_2=v^2_1$. It has the form
\begin{equation}\label{7d}
\frac{r}{r+1}\varphi_2=\varphi'_1.
\end{equation}
We take the particular solution of reduced system of ordinary differential equations (\ref{6d}), (\ref{7d}) in the form
\[
\varphi_1=\left (\sqrt{\frac{r(r+1)}{2(r-1)}}\omega+C_1 \right )^{\frac{2}{r+1}},
\]
\[
\varphi_2=\sqrt{\frac{2(r+1)}{r(r-1)}}\left (\sqrt{\frac{r(r+1)}{2(r-1)}}\omega+C_1 \right )^{\frac{1-r}{r+1}}
\]
where $C_1=const$. Thus one has to integrate overdetermined compatible system of differential equations
\[
u_{x_1}=(x_1)^{\frac{-1}{r+1}}\left (\sqrt{\frac{r(r+1)}{2(r-1)}}(x_2-u_{x_2})+C_1 \right )^{\frac{2}{r+1}}
\]
\[
u_{x_2}=(x_1)^{\frac{r}{r+1}}\sqrt{\frac{2(r+1)}{r(r-1)}}\left (\sqrt{\frac{r(r+1)}{2(r-1)}}(x_2-u_{x_2})+C_1 \right )^{\frac{1-r}{r+1}}
\]
to construct the solution of equation (\ref{2}).
Nevertheless it is easy to get the solution of the equation
\begin{equation}\label{8d}
w_t- \frac{1}{r}\left (\frac{(w-1)^{\frac{1-r}{r}}w_x}{w^{\frac{1+r}{r}}}\right )_x=0
\end{equation}
in such form
\[
\left (1-\frac{1}{w} \right )^{\frac{1}{r}}=(t)^{\frac{-1}{r+1}}\left (\sqrt{\frac{r(r+1)}{2(r-1)}}(x-\theta )+C_1 \right )^{\frac{2}{r+1}},
\]
\[
\theta=(t)^{\frac{r}{r+1}}\sqrt{\frac{2(r+1)}{r(r-1)}}\left (\sqrt{\frac{r(r+1)}{2(r-1)}}(x-\theta)+C_1 \right )^{\frac{1-r}{r+1}}.
\]

The tangent transformations groups are also used in the framework of this approach. Let us consider the nonlinear evolution equation
\begin{equation}\label{9d}
u_t=e^{\frac{1}{u_{xx}}}.
\end{equation}
One can construct operator of tangent transformations of the form
\[
K= -t\partial_{t}+u_x\partial_{x}+\frac{u_x^2}{2}\partial_u+ u_t\partial_{u_t}
\]
admitted by (\ref{9d}). The first order functionally independent differential invariants of the corresponding one-parameter Lie group of tangent transformations can be chosen in the form
\[
\omega=xu_x-2u, \quad \omega_1=\ln u_t-\frac{x}{u_x}, \quad \omega_2=u_x, \quad \omega_3=tu_t.
\]
In order to construct ansatz of type (\ref{1}) reducing equation (\ref{9d}) to system of ordinary differential equations we consider two-dimensional Lie algebra with basic operators $\{K, P_t=\frac{\partial}{\partial t} \}$. The operators satisfy the commutation relation $[K,P_t]=P_t$. The invariants of two-parameter Lie group with generators $K, P_t$ are $\omega, \omega_1, \omega_2$. Then  we construct the ansatz by using these invariants in the form
\begin{equation}\label{11d}
u_x=f(\omega), \quad u_t=\exp \left (\varphi(\omega)+\frac{x}{f(\omega)}\right ).
\end{equation}
From (\ref{11d}) and (\ref{9d}) we have
\[
u_{xx}=-\frac{f'f}{1-xf'}
\]
and first ordinary differential equation
\begin{equation}\label{12d}
f'f\varphi =-1.
\end{equation}
From the condition $u_{xt}=u_{tx}$ it follows that $f, \varphi$ satisfy the second ordinary differential equation
\begin{equation}\label{13d}
\frac{f-\varphi' f^3}{f^2}=-2f'.
\end{equation}
Thus the reduced system consists of equations (\ref{12d}) and (\ref{13d}). From (\ref{12d}), (\ref{13d}) and (\ref{11d}) it follows that the solutions of equation (\ref{9d}) can be constructed by integrating overdetermined compatible system
\[
u_t=\exp{\left [\frac{2(\sqrt{C_1-4(xu_x-2u)}+x)}{2C_2-\sqrt{C_1-4(xu_x-2u)}}\right ]},
\]
\[
u_x=\frac{2C_2-\sqrt{C_1-4(xu_x-2u)}}{2}
\]
where $C_1$, $C_2$ are arbitrary real constants.

Next we emphasize  that the operators of conditional symmetry of corresponding system can be used for construction the B\"acklund transformations for nonlinear wave equation
\begin{equation}\label{14d}
u_{x_1x_2}=[1-k^2u_{x_2}^2]^{1/2}\sin u.
\end{equation}
Indeed we showed that
\begin{equation}\label{14d1}
Q=\partial_{x_3}+k\cos x_3\partial_{v^1}+k^{-1}\sqrt{1-k^2(v^2)^2}\partial_{v^2}
\end{equation}
is the operator of conditional symmetry of the corresponding system
\begin{equation}\label{8}
v^1_2+v^1_3v^2=v^2_1+v^2_3v^1,
\end{equation}
\begin{equation}\label{9}
v^2_1+v^2_3v^1=\sqrt{1-k^2(v_2)^2}\sin x_3,
\end{equation}
where $u\equiv x_3$. Using operator $Q$ we can write the ansatz in the following form
\begin{equation}\label{15d}
u_{x_1}=\varphi_2+k\sin u, \quad u_{x_2}=k^{-1}\sin (u-\varphi_1),
\end{equation}
where $\varphi_1$, $\varphi_2$ are unknown functions  on $x_1$, $x_2$ and hence the B\"acklund transforms
\begin{equation}\label{16d}
u_{x_2}=k^{-1}\sin (u-w), \quad u_{x_1}=w_{x_1}+k\sin u
\end{equation}
relating equation (\ref{14d}) and sine-Gordon equation $w_{x_1x_2}=\sin w$.
These B\"acklund transforms (\ref{16d})  have been obtained for the first time in \cite{BD} by another technique.

Note that this approach is also applicable for  linearization nonlinear partial differential equations with two independent variables. Indeed, consider the second order differential equation
\begin{equation}\label{17d}
u_{x_0x_0}=F(u_{x_0x_1,} u_{x_1x_1}),
\end{equation}
where $F$ is a smooth function. Using the invariance  of (\ref{17d}) under Lie group of transformations with corresponding five-dimensional Lie algebra given by basic
elements $\partial_{x_0}$, $\partial_{x_1}$ $\partial_u$, $x_0\partial_u$, $x_1\partial_u$ we write the corresponding system in the form
\begin{equation}\label{18d}
\frac{\partial F}{\partial v^2} \frac{\partial v^2}{\partial x_1}+\frac{\partial F}{\partial v^3} \frac{\partial v^3}{\partial x_1} =\frac{\partial v^2}{\partial x_0}, \quad \frac{\partial v^3}{\partial x_0}=\frac{\partial v^2}{\partial x_1},\quad v^1=F(v^2, v^3),
\end{equation}
where $u_{x_0x_0}\equiv v^1(x_0,x_1)$, $u_{x_0x_1}\equiv v^2(x_0,x_1)$, $u_{x_1x_1}\equiv v^3(x_0,x_1)$. One can prove that (\ref{18d}) possesses infinite Lie classical symmetry and can be linearized by hodograph transformations.
Thus we obtained the method of linearization of the second order partial differential equation of the form (\ref{17d}) for arbitrary function $F$.

Let us note that the symmetry group of corresponding system written in the general form contains the symmetry group of point transformations of initial equation as a subgroup and generators of point transformations can be used to construct ansatz (\ref{1}). However these operators lead to invariant solutions in the classical Lie sense. We shall illustrate this property by the following example. Let us consider the wave equation
\begin{equation}\label{19d}
u_{x_1x_2}=F(u),
\end{equation}
where $F$ is a smooth function. It is invariant with respect to the three-parameter Lie group. The basis of Lie algebra is  given by
$\{\partial_{x_1}, \partial_{x_2}, x_1\partial_{x_1}-x_2\partial_{x_2} \}$. Consider two-dimensional subalgebra with basic elements $\{ \partial_{x_2}, x_1\partial_{x_1}-x_2\partial_{x_2} \}$.
 By using the differential invariants   $u, x_1u_{x_1}, \frac{u_{x_2}}{x_1}$ of the corresponding two-parameter Lie group we construct ansatz of the form
\begin{equation}\label{20d}
u_{x_1}=\frac{f(u)}{x_1}, \quad u_{x_2}=x_1\varphi (u)
\end{equation}
which reduces (\ref{19d}) to the system
\[
f'\varphi=\varphi +\varphi' f=F(u).
\]
Let $F(u)=0$. Then we obtain two cases\\
\[
1. f'=0, \quad \varphi +\varphi' f=0
\]
and solution of reduced system has the form
\[
f=C_1=const, \quad \varphi= C_2\exp{\left (-\frac{u}{C_1}\right )}, \quad C_2=const.
\]
By integrating system
\[
\frac{\partial u}{\partial x_1}=\frac{C_1}{x_1}, \quad \frac{\partial u}{\partial x_2}=C_2x_1\exp{\left (-\frac{u}{C_1} \right )}
\]
one obtains the solution
\begin{equation}\label{21d}
u=C_1\ln \left (\frac{C_2}{C_1}x_1x_2+C_3x_1\right )
\end{equation}
where $C_3$ is arbitrary real constant and $C_1\not =0$, of equation (\ref{19d}) with $F=0$. In the second case we have
\[
2. \hspace{5mm} \varphi =0, 
\]
\[
u_{x_1}=\frac{1}{x_1}f(u), \quad u_{x_2}=0
\]
and solution has the form
\begin{equation}\label{22d}
u=h(x_1),
\end{equation}
where $h(x_1)$ is arbitrary differentiable function.
Let consider the operator
\begin{equation}\label{23d}
Q=\alpha \partial_{x_2}+\beta ( x_1\partial_{x_1}-x_2\partial_{x_2})
\end{equation}
where $\alpha$, $\beta$ are arbitrary real constants. One can verify that
\[
Q\left (u-C_1\ln \left (\frac{C_2}{C_1}x_1x_2+C_3x_1\right )\right )=0
\]
iff
\begin{equation}\label{24d}
\alpha \frac{C_2}{C_1}+\beta C_3=0.
\end{equation}
 It means  that solution (\ref{21d}) is  invariant  with respect to one-dimensional subgroup of symmetry group of equation (\ref{19d}) with generator $Q$ where $\alpha$, $\beta$ satisfy the condition (\ref{24d}). It is obvious that the solution (\ref{22d}) is  invariant  with respect to one-parameter group with generator  $Q=\alpha \partial_{x_2}$ ($\beta=0$). Thus we conclude that any solution of equation (\ref{19d}) when $F=0$ constructed by this method with the help of two-dimensional Lie algebra with basic elements $\{ \partial_{x_2}, x_1\partial_{x_1}-x_2\partial_{x_2} \}$ is an invariant one in the classical Lie sense. 

Further we show how the operators of Lie--B\"acklund symmetry \cite{O,Fok} are used for reducing partial differential equations. Let us consider equation
\begin{equation}\label{25d}
U ( x, u,{\mathop u\limits_1},{\mathop u\limits_2},\ldots ,
{\mathop u\limits_k} ) =0,
\end{equation}
where $x=(x_1,x_2,\ldots ,x_n)$, $u=u(x)\in C^k(\mathbb{R}^n,{\mathbb
R}^1)$, and
${\mathop u\limits_k}$ denotes all partial derivatives of $k$-th order and the $m$-th order ordinary differential equation of the form
\begin{equation}\label{26d}
H\left (x_1, x_2, \ldots , x_n, u, \frac{\partial u}{\partial x_1}, \ldots, \frac{\partial^m u}{\partial x_1^m} \right )=0.
\end{equation}
Let
\[
u=F(x, C_1,\ldots,C_{m}),
\]
where $F$ is a smooth function on variables $x, C_1,\ldots ,C_{m}$, and
$C_1,\ldots, C_{m}$
are arbitrary functions on variables $x_2, x_3,\ldots, x_n$,
be a general solution of equation~(\ref{26d}).
 
We use the Theorem 1 from \cite{T8} which implies that if 
 equation~(\ref{26d}) is invariant with respect to
the Lie--B\"acklund operator $X=U ( x, u,{\mathop u\limits_1},{\mathop u\limits_2},\ldots ,{\mathop u\limits_k} )\partial_u$ then the  ansatz
\begin{equation}\label{28d}
u=F(x, \varphi_1, \varphi_2,\ldots,\varphi_{m}),
\end{equation}
where $\varphi_1, \varphi_2,\ldots,\varphi_{m}$ depend on $n-1$ variables
 $x_2, x_3,\ldots, x_n$ reduces partial differential equation (\ref{25d})
to the system of $k_1$ equations for unknown functions
 $\varphi_1, \varphi_2,\ldots ,\varphi_m$ with $n-1$ independent variables
 and $k_1 \le m$.
We show  the application of the theorem to nonlinear partial differential equation.
Consider linear ordinary differential equation
\begin{equation}\label{29d}
u_{x_1x_1}+\alpha^2u_{x_1}=0
\end{equation}
where $\alpha=const$. Recall, that the concepts of local theory of differential equations such as symmetry, conditional symmetry, conservation laws, Lax representations are defined by differential equalities which must be satisfied only for solutions of the equations under study. One can prove that equation (\ref{29d}) admits the following Lie--B\"acklund operator
\[
X=(u_{x_1x_2}-u_{x_1}F(u_{x_1}+\alpha^2u))\partial_u,
\]
where $F\in C^2(\mathbb{R}^1,{\mathbb
R}^1)$. It means that the  following criterium of invariance
\[
X_{(2)}(u_{x_1x_1}+\alpha^2u_{x_1})=0 \quad \textrm{whenever} \quad u_{x_1x_1}+\alpha^2u_{x_1}=0,
\]
where $X_{(2)}$ is the prolongated operator of the second order \cite{O}, is fulfilled. 
We have proved that (\ref{29d}) admits operator of non-point (tangent) symmetry 
\[
X_1=f(u, u_x)\partial u
\]
if $f(u, u_x)$ satisfies the following equation
\[
f''_{uu} -2\alpha^2f''_{uu_x}+ \alpha^4f''_{u_xu_x}=0.
\]
The general solution of this equation has the form
\begin{equation}\label{1m}
f=A(u_{x_1}+\alpha^2u)u+B(u_{x_1}+\alpha^2u)
\end{equation}
where $A$, $B$ are arbitrary smooth functions of one variable. One can verify that equation (\ref{29d}) also admits  operator 
\[
X_2=e^{-\alpha^2x_1}h(u_{x_1}+\alpha^2u)\partial u
\]
where $h$ is arbitrary function on variable $u_{x_1}+\alpha^2u$.
 Then the ansatz
\begin{equation}\label{30d}
u=\varphi_1(x_2)+e^{-\alpha^2x_1}\varphi_2(x_2)
\end{equation}
obtained from the general solution of equation (\ref{29d})
reduces wave type partial differential equations
\begin{equation}\label{31d}
u_{x_1x_2}=u_{x_1}F(u_{x_1}+\alpha^2u)+A(u_{x_1}+\alpha^2u)u+B(u_{x_1}+\alpha^2u)+ ku_{x_2}+e^{-\alpha^2x_1}h(u_{x_1}+\alpha^2u)
\end{equation}
where  $k$ is a real constant. In general, the $x_1$ dependent coefficients in partial differential equations enable us to study the effects of field gradients.

Substituting (\ref{30d}) into (\ref{31d}) we obtain the reduced system of two ordinary differential equations
\begin{equation}\label{32d}
-\alpha^2\varphi'_2=-\alpha^2\varphi_2F(\alpha^2\varphi_1)+A(\alpha^2\varphi_1)\varphi_2+k\varphi'_2+h(\alpha^2\varphi_1),
\end{equation}
\begin{equation}\label{32m}
A(\alpha^2\varphi_1)\varphi_1+ B(\alpha^2\varphi_1)+k\varphi'_1=0 
\end{equation}
for unknown functions $\varphi_1(x_2)$, $\varphi_2(x_2)$. One can obtain partial solutions of (\ref{31d}) from solutions of system (\ref{32d}), (\ref{32m}). In particular, if $A=B=0$ and $k=0$ then system (\ref{32d}), (\ref{32m})
is reduced to one ordinary differential equation of the form
%\begin{equation}\label{33m}
\[
\varphi'_2=\varphi_2F(\alpha^2\varphi_1)-\displaystyle \frac{1}{\alpha^2}h(\alpha^2\varphi_1) 
\]
%\end{equation}
where $\varphi_1(x_2)$, $\varphi_2(x_2)$ are unknown functions.
This equation is integrable by quadratures for arbitrary $\varphi_1(x_2)$. Its general solution has the form
\begin{equation}\label{35m}
\varphi_2=\left (C_1-\displaystyle \frac{1}{\alpha^2}\int h(\alpha^2\varphi_1(x_2))H(x_2)dx_2 \right )
\exp{\left (\int F(\alpha^2\varphi_1(x_2))dx_2 \right )}
\end{equation}
where $C_1=const$,
\begin{equation}\label{35v}
H(x_2)=\exp{\left (-\int F(\alpha^2\varphi_1(x_2))dx_2 \right )}.
\end{equation}
Using (\ref{30d}) one can construct the solution of nonlinear wave equation 
\begin{equation}\label{31w}
u_{x_1x_2}=u_{x_1}F(u_{x_1}+\alpha^2u)+e^{-\alpha^2x_1}h(u_{x_1}+\alpha^2u)
\end{equation}

 in the following form
\begin{equation}\label{33d}
u=\varphi_1(x_2)+\left (C_1-\displaystyle \frac{1}{\alpha^2}\int h(\alpha^2\varphi_1(x_2))H(x_2)dx_2 \right )
\exp{\left (\int F(\alpha^2\varphi_1(x_2))dx_2 -\alpha^2 x_1\right )},
\end{equation}
where $\varphi_1(x_2)$ is arbitrary smooth function. So, in the framework of this approach we have constructed
solution with arbitrary function $\varphi_1(x_2)$ to nonlinear wave type partial differential equation (\ref{31w})
for arbitrary functions $F$ and $h$.

\section{Conclusions}
We have constructed  ans\"atze (\ref{4d}) and (\ref{11d}) which reduce nonlinear evolution equations (\ref{2}) and (\ref{9d}) to ordinary differential equations and can not be obtained by using classical Lie method. We have found  the solution of nonlinear heat equation (\ref{8d}). It turns out that some of these ans\"atze result in the classical invariant solutions.  Obviously, one can  construct such ans\"atze  by prolongated operators of point symmetry admitted by the initial equation but they lead to the  invariant solutions too. It is necessary that operators of non-point and conditional symmetry should be applied to obtain new results.

As was noted above the linearization of class of nonlinear partial differential equations (\ref{17d}) is possible in the framework of this approach.

Finally we show that the existence even at least one operator of Lie-B\"acklund symmetry to ordinary differential equations (\ref{29d})  gives the possibility of constructing the solutions  (\ref{33d})  defined by arbitrary functions to equation (\ref{31w}). To our knowledge the inverse scattering tranformation method is not applicable in this case.
We show that this approach is applicable to nonevolutionary differential equations.


\begin{thebibliography}{99}

\small

\bibitem{BK}
Bluman G.W. and Kumei S.,
Symmetries and Differential Equations, (1989), Appl. Math. Sci. 81, Springer--Verlag, Berlin

\bibitem{O}
Olver P.J., Applications of Lie Groups to Differential
Equations, Springer--Verlag, New York, 1986

\bibitem{Ovs}
Ovsiannikov L.V., Group Analysis of Differential Equations, Academic Press, New York, 1982 

\bibitem{Fok} 
Fokas A.S. and Q.M.Liu, Nonlinear interaction of travelling waves of nonintegrable equations, {\it Phys. Rev. Lett.} {\bf 72, No.21} (1994), 3293--3296.
\bibitem{Zhd}
Zhdanov R.Z., 
Conditional Lie-B\"acklund symmetry and reduction of evolution equations,
{\it J. Phys. A: Math. Gen. } {\bf 28} (1995), 3841--3850.

\bibitem{Kun}
 Kunzinger M.  and  Popovych R., Generalized conditional symmetries of evolution equations,
{\it J. Math. Anal. Appl} {\bf 379} (2011), 444-460.

\bibitem{BluKyt}
Bluman G.W. and Zhengzheng Yang, A symmetry-based method for constructing nonlocally related partial differential equation systems, {\it J. Math. Phys} {\bf 54}, 093504 (2013).

\bibitem{T2}
Tsyfra I., Napoli A., Messina A., Tretynyk V., On new ways of group methods for reduction of evolution-type
equations ,
{\it J. Math. Anal. Appl.} {\bf 307}, (2005), 724--735.

\bibitem{BD}
Dodd R.K. and Bullogh R.K. B\"acklund transformations for the sine--Gordon equations, {\it Proc. R. Soc. A} {\bf 351}, (1976), 499--523.

\bibitem{T8} Tsyfra I.M. Symmetry reduction of nonlinear differential equations, Proceedings of Institute of Mathematics, Kiev, 2004, {\bf 50}, pp. 266--270.



\end{thebibliography}
\end{document}